\documentclass{article}

\usepackage{mystyle}
\usepackage{wasysym}
\usepackage{mathtools}
\begin{document}
\title{A note on Woodin's HOD dichotomy}
\author{Gabriel Goldberg\\ Evans Hall\\ University Drive \\ Berkeley, CA 94720}
\maketitle
\section{Introduction}
The purpose of this note is to prove a version of Woodin's HOD dichotomy
\cite{Woodin}
from a strongly compact cardinal.
A cardinal \(\kappa\) is strongly compact
if every \(\kappa\)-complete filter extends to a \(\kappa\)-complete
ultrafilter. The class of hereditarily ordinal definable sets, denoted \(\HOD\),
is the minimum class \(M\) that contains all the ordinals 
and is definably closed:
if \(S\subseteq M\) is definable (in the universe of sets \(V\)) by a set theoretic formula
with parameters in \(M\), then \(S\in M\). \(\HOD\) is a proper class transitive model
of ZFC.
\begin{thm}\label{thm:main}
    Suppose \(\kappa\) is strongly compact. 
    Then one of the following holds:
    \begin{enumerate}[(1)]
        \item For every singular strong limit cardinal \(\lambda\geq \kappa\),
        \(\lambda\) is singular in \(\HOD\) and \(\lambda^{+\HOD} = \lambda^+\).\label{item:closetov}
        \item All sufficiently large regular cardinals are 
        measurable in \(\HOD\).\label{item:farfromv}
    \end{enumerate}
\end{thm}

If \ref{item:closetov} holds, \ref{item:farfromv} cannot (since there are arbitrarily large
successor cardinals that are not inaccessible in \(\HOD\)), so \cref{thm:main} is
truly a dichotomy. 

We will actually prove a stronger dichotomy. An inner model \(M\)
is said to have the \(\lambda\)-cover property if every set of ordinals
of cardinality less than \(\lambda\) is covered by (that is, included in)
a set of ordinals in \(M\) of cardinality less than \(\lambda\).

\begin{thm}\label{thm:Jensen}
    Suppose \(\kappa\) is strongly compact. Then one of the following holds:
    \begin{enumerate}[(1)]
        \item For every strong limit cardinal \(\lambda \geq \kappa\), 
        \(\HOD\) has the \(\lambda\)-cover property.\label{item:closetov1}
        \item All sufficiently large regular cardinals are \(\omega\)-strongly 
        measurable in \(\HOD\).\label{item:farfromv1}
    \end{enumerate}
\end{thm}
The definition
of \(\omega\)-strong measurability is deferred until the beginning
of the next section, but let us note here that any cardinal 
that is \(\omega\)-strongly measurable in \(\HOD\) is
indeed measurable in \(\HOD\). Note that
\cref{thm:Jensen} \ref{item:closetov1} strengthens 
\cref{thm:main}  \ref{item:closetov}, while
\cref{thm:Jensen} \ref{item:farfromv1} strengthens 
\cref{thm:main}  \ref{item:farfromv}. 

Assuming instead that \(\kappa\) is \(\HOD\)-supercompact,
Woodin \cite{Woodin} proved an even stronger dichotomy theorem,
for example establishing that if (either) condition \ref{item:closetov} above holds,
then \(\kappa\) is supercompact in \(\HOD\).
On the other hand, Cheng-Friedman-Hamkins \cite{ChengFriedmanHamkins}
produce a model of ZFC with a supercompact cardinal \(\kappa\)
in which \ref{item:closetov} holds, yet 
\(\kappa\) is not weakly compact in \(\HOD\), which shows that
Woodin's theorem cannot be proved from the hypotheses we assume here.
\section{The proof}
Suppose \(\delta\) is a regular cardinal and \(S\subseteq \delta\)
is a stationary set. We say an inner model \(M\) splits \(S\)
if \(S\in M\) and for all \(\gamma\) such that \((2^\gamma)^M < \delta\),
there is a partition of \(S\) into \(\gamma\)-many stationary sets.
A cardinal is \(\omega\)-strongly measurable
in \(\HOD\) if \(\HOD\) does not split the set
of ordinals less than \(\delta\) 
with countable cofinality, which is denoted by \(S^\delta_\omega\).
\begin{thm}[{Woodin, \cite[Lemma 3.37]{Woodin}}]\label{thm:Woodin}
    If \(S\) is an ordinal definable stationary subset of
    a regular cardinal \(\delta\) and \(\HOD\) does not split \(S\),
    then \((\mathcal C_\delta\restriction S)\cap \HOD\) is atomic in \(\HOD\).
    In particular, if \(\delta\) is \(\omega\)-strongly measurable in \(\HOD\),
    then \(\delta\) is measurable in \(\HOD\) and
    \(\delta\) contains an \(\omega\)-club of cardinals inaccessible
    in \(\HOD\). 
\end{thm}
The last sentence is supposed to highlight
that the existence of \(\omega\)-strongly measurable cardinals
entails a massive failure of the cover property
even for countable sets. 
Actually there is an \(\omega_1\)-club of cardinals
\textit{measurable} in \(\HOD\)
below any cardinal that is \(\omega_1\)-strongly measurable
in \(\HOD\) in the natural sense.

\begin{prp}\label{prp:concentration}
    Suppose \(\kappa\) is strongly compact,
    \(\delta\geq \kappa\) is regular, and \(\HOD\) splits 
    \(S^\delta_\omega\). Then for any \(\gamma\) such that
    \((2^\gamma)^\HOD < \delta\), there is a \(\kappa\)-complete
    fine ultrafilter \(\mathcal U\) on \(P_\kappa(\gamma)\)
    such that \(\mathcal U\) concentrates on \(P_\kappa(\gamma)\cap \HOD\).
    \begin{proof}
        Let \(j : V\to M\) be an elementary embedding from the universe into an
        inner model such that \(\crit(j) = \kappa\)
        and \(j[\delta]\) is contained in some
        set \(S\in M\) with \(|S|^M < j(\kappa)\). In particular,
        the ordinal \(\delta_* = \sup j[\delta]\) has cofinality
        less than \(j(\kappa)\) in \(M\). Let \(C\subseteq \delta_*\)
        be a closed unbounded set of ordertype \(\cf^M(\delta_*)\).

        Fix a cardinal \(\gamma\) such that
        \((2^\gamma)^\HOD < \delta\)
        and let \(\langle S_\alpha \rangle_{\alpha < \gamma}\)
        witness that \(\HOD\) splits \(S^\delta_\omega\).
        Then let 
        \[\sigma = \{\xi < j(\gamma): T_\xi\cap \delta_*\text{ is stationary in }M\}\]
        where \(\vec T = j(\vec S)\).
        Thus \(\sigma\in \HOD^M\). Notice that \(j[\gamma]\subseteq \sigma\):
        \(j[S_\xi]\subseteq T_{j(\xi)}\), so in fact 
        \(T_{j(\xi)}\cap \delta_*\) is truly stationary
        (not just in \(M\)). For all \(\xi \in \sigma\), 
        \(T_\xi\cap C\neq \emptyset\), so let \(f(\xi) = \min(T_\xi\cap C)\).
        Then \(f :\sigma  \to C\) is an injection.
        So \(|\sigma| = \cf^M(\delta_*)\). In particular,
        \(\sigma\in j(P_\kappa(\gamma))\).

        Finally let \(\mathcal U\) be the ultrafilter
        on \(P_\kappa(\gamma)\) derived from \(j\) using \(\sigma\).
        That is, let \(\mathcal U = \{A\subseteq P_\kappa(\gamma) : \sigma\in j(A)\}\).
        Since \(\sigma\in \HOD^M\), \(\sigma \in j(P_\kappa(\gamma)\cap \HOD)\),
        and hence \(\HOD\cap P_\kappa(\gamma)\in \mathcal U\).
        Since \(j[\gamma]\subseteq \sigma\), \(\mathcal U\) is fine.
        Since \(\crit(j) = \kappa\), \(\mathcal U\) is \(\kappa\)-complete.
    \end{proof}
\end{prp}
The main observation involved in the proof above is that
the stationary splitting argument from \cite{Woodin}
(which Woodin calls ``Solovay's Lemma" although it is a bit
different from the related lemma in \cite{Kanamori}) can be
adapted to strongly compact cardinals. Usuba \cite{USUBA2020} made the
same observation independently and earlier.

\begin{lma}\label{KappaCov}
    Suppose \(\kappa\) is strongly compact. 
    Then one of the following holds:
    \begin{enumerate}[(1)]
        \item \(\HOD\) has the \(\kappa\)-cover property.\label{item:closetov}
        \item All sufficiently large regular cardinals are 
        \(\omega\)-strongly measurable in in \(\HOD\).\label{item:farfromv}
    \end{enumerate}
\begin{proof}
    Assume that there are arbitrarily large regular cardinals \(\delta\) that are not \(\omega\)-strongly 
    measurable in \(\HOD\), or in other words,
    \(\HOD\) splits \(S^\delta_\omega\).
    Applying \cref{prp:concentration} to sufficiently large such \(\delta\),
    for all \(\gamma\geq \kappa\), there is a \(\kappa\)-complete
    fine ultrafilter \(\mathcal U\) on \(P_\kappa(\gamma)\) such that
    \(P_\kappa(\gamma)\cap \HOD\in \mathcal U\). For each \(\sigma\in P_\kappa(\gamma)\),
    \(\{\tau\in P_\kappa(\gamma) :\sigma\subseteq \tau\}\in \mathcal U\). Since \(\mathcal U\)
    is a filter, it follows that 
    \(\{\tau\in P_\kappa(\gamma)\cap \HOD :\sigma\subseteq \tau\}\in \mathcal U\),
    and in particular this set is nonempty. This yields \(\tau\in P_\kappa(\gamma)\cap \HOD\)
    covering \(\sigma\), as desired.
\end{proof}
\end{lma}
We now extend the cover property of \(\HOD\) to all strong limit cardinals greater than or equal to
the first strongly compact cardinal.
\begin{proof}[Proof of \cref{thm:Jensen}]
        Suppose \ref{item:farfromv} fails, 
        so by \cref{KappaCov}, \(\HOD\) has the \(\kappa\)-cover property.
        Given this, it suffices to show that for all \(\delta \geq \kappa\), 
        for some \(\delta' < \beth_\omega(\delta)\), every set of ordinals of cardinality
        at most \(\delta\) is covered by a set of ordinals in \(\HOD\) of cardinality at
        most \(\delta'\).

        Let \(\mathcal U\) be a fine \(\kappa\)-complete ultrafilter on \(P_\kappa(\delta)\).
        Fix \(A\subseteq \kappa\) such that \(V_\kappa\subseteq \HOD_A\).
        We will show that for every set \(S\) of ordinals of cardinality at most \(\delta\),
        there is a set \(T\) of cardinality at most \(2^\delta\) such that \(T\in \HOD_{A,\mathcal U}\)
        and \(S\subseteq T\). By Vopenka's theorem, \(\HOD_{A,\mathcal U}\) is a forcing extension
        of \(\HOD\) by a forcing in \(\HOD\) of cardinality less than \(\beth_\omega(\delta)\),
        so the theorem follows.

        Note that \(\HOD_{A,\mathcal U}\) is closed under \({<}\kappa\)-sequences
        since \(\HOD_{A,\mathcal U}\) has the \(\kappa\)-cover property and 
        \(V_\kappa \subseteq \HOD_{A,\mathcal U}\).
        As a consequence, \(M_\mathcal U\) satisfies that \(j_\mathcal U(\HOD_{A,\mathcal U})\) 
        is closed under \({<}j_\mathcal U(\kappa)\)-sequences.
        Also \(j_\mathcal U(\HOD_{A,\mathcal U})\subseteq \HOD_{A,\mathcal U}\): there is a wellorder of
        \(j_\mathcal U(\HOD_{A,\mathcal U})\) definable from \(A\) and \(\mathcal U\).
        
        Finally, suppose \(S\) is a set of ordinals of cardinality at most \(\delta\). Then 
        since \(\mathcal U\) is fine, there is some
        \(T\) in \(M_\mathcal U\) such that \(S\subseteq T\) and \(|T|^{M_\mathcal U} < j_{\mathcal U}(\delta)\). 
        Therefore \(T\in j_\mathcal U(\HOD_{A,\mathcal U})\subseteq \HOD_{A,\mathcal U}\). 
        Moreover \(|T| \leq |j_\mathcal U(\delta)| \leq 2^\delta\). This completes the proof.
\end{proof}

It may be unclear where we used that \(\kappa \neq \omega\) in the proof above.
In fact, we used the countable completeness of \(\mathcal U\) to establish that
there is a wellorder of \(j_\mathcal U(\HOD_{A,\mathcal U})\) definable from \(A\) and \(\mathcal U\).
(If \(\mathcal U\) is countably incomplete, the canonical wellorder of 
\(j_\mathcal U(\HOD_{A,\mathcal U})\) as computed in \(M_\mathcal U\)
is not a wellorder at all.)

Using these theorems, we prove the dichotomy involving successors of singulars.
\begin{proof}[Proof of \cref{thm:main}]
    By \cref{KappaCov} (and \cref{thm:Woodin}), 
    we can assume that \(\HOD\) has the \(\kappa\)-cover property.
    Fix a singular strong limit cardinal \(\lambda > \kappa\). 
    \cref{thm:Jensen} easily implies that \(\lambda\) is singular in \(\HOD\). 
    Let \(\gamma = \lambda^{+\HOD}\). 
    Since \(\gamma\) is regular in \(\HOD\),
    \cref{thm:Jensen} implies \(\cf (\gamma) \geq \lambda\).
    Since \(\lambda\) is singular, \(\cf(\gamma) > \lambda\),
    and so \(\gamma = \lambda^+\).
\end{proof}

\bibliographystyle{plain}
\bibliography{Bibliography.bib}

\begin{thebibliography}{1}

\bibitem{ChengFriedmanHamkins}
Yong Cheng, Sy-David Friedman, and Joel~David Hamkins.
\newblock Large cardinals need not be large in {HOD}.
\newblock {\em Annals of Pure and Applied Logic}, 166(11):1186 -- 1198, 2015.

\bibitem{Kanamori}
Robert~M. Solovay, William~N. Reinhardt, and Akihiro Kanamori.
\newblock Strong axioms of infinity and elementary embeddings.
\newblock {\em Ann. Math. Logic}, 13(1):73--116, 1978.

\bibitem{USUBA2020}
Toshimichi Usuba.
\newblock A note on $\delta$-strongly compact cardinals.
\newblock {\em Topology and its Applications}, page 107538, 2020.

\bibitem{Woodin}
W.~Hugh Woodin.
\newblock In search of {U}ltimate-{$L$}: the 19th {M}idrasha {M}athematicae
  {L}ectures.
\newblock {\em Bull. Symb. Log.}, 23(1):1--109, 2017.

\end{thebibliography}
\end{document}